\documentstyle[amssymb,amsfonts, 11 pt]{amsart}

\newtheorem{thm}{Theorem}[section]
\newtheorem{cor}[thm]{Corollary}
\newtheorem{lem}[thm]{Lemma}
\newtheorem{prop}[thm]{Proposition}
\theoremstyle{definition}
\newcommand{\comment}[1]{}

\theoremstyle{remark}

\numberwithin{equation}{section}

\def\be#1{\begin{equation}\label{#1}}

\newenvironment{proof}{\noindent {\bf Proof} }{\endprf\par}
\def \endprf{\hfill  {\vrule height6pt width6pt depth0pt}\medskip}
\def\emph#1{{\it #1}}
\def\textbf#1{{\bf #1}}
\begin{document}       

\title[Navier-Stokes in $BMO^{-1}$]{Regularity of solutions to the
  Navier-Stokes equations evolving from small data in $BMO^{-1}$}

\author[P. Germain]{Pierre Germain} 
\address{Courant Institute of Mathematical Sciences, New York University, New York, NY 10012-1185}  
\email{\tt pgermain@@math.nyu.edu}

\author[N. Pavlovi\'{c}]{Nata\v{s}a Pavlovi\'{c}}
\address{Department of Mathematics, Princeton University, Princeton, NJ 08544-1000}
\email{\tt natasa@@math.princeton.edu}

\author[G. Staffilani]{Gigliola Staffilani}
\address{Department of Mathematics, Massachusetts Institute of Technology,
Cambridge, MA 02139-4307}
\email{\tt gigliola@@math.mit.edu}
\thanks{N.P. was supported in part by N.S.F. Grant DMS 0304594 and
G.S. was supported in part by N.S.F. Grant DMS 0602678.}

\date{January 23, 2007}

\begin{abstract}
In 2001, H. Koch and D. Tataru proved the existence of
global in time solutions to the incompressible Navier-Stokes equations in
${\mathbb{R}}^d$ for initial data small enough in $BMO^{-1}$. 
We show in this article that the Koch and Tataru
solution has higher regularity. As a consequence, we get a decay estimate 
in time for any space derivative, and space analyticity of the solution.
Also as an application of our regularity theorem, we prove a regularity result
for self-similar solutions. 
\end{abstract}

\maketitle

\section{Introduction}

In this paper we study regularity, decay and analyticity 
of solutions to the Navier-Stokes equations 
for the incompressible fluid in ${\mathbb R}^{d}$, which are given by
\be{ns} \frac{\partial u}{\partial t} + (u \cdot \nabla) u + \nabla p
= \Delta u + f, 
\end{equation}
\be{nsdiv} \nabla \cdot u = 0,
\end{equation}
and the initial condition
\be{0nsin} u(x,0) = u_{0}(x),
\end{equation}
for the unknown velocity vector field $u = u(x,t) \in {\mathbb R}^{d}$ and 
the pressure $p = p(x,t) \in {\Bbb R}$, where $x \in {\mathbb R}^{d}$ and 
$t \in [0,\infty)$. In the rest of this note we shall take $f(x,t) = 0$. 

Existence of global in time solutions to \eqref{ns}-\eqref{0nsin}
when the space dimension $d = 3$, their uniqueness and regularity 
are long standing open problems of fluid dynamics. An approach in 
studying solutions to the Navier-Stokes equations is to construct 
solutions to the corresponding integral equation via a fixed point theorem, 
so called ``mild'' solutions.  In the context of the
Navier-Stokes equation this approach was pioneered by Kato and Fujita, see for example, 
\cite{KF}. However the existence of mild solutions to the Navier-Stokes equations
\eqref{ns} - \eqref{0nsin} in ${\mathbb R}^{d}$ for $d \geq 3$ has been established only locally in time
and globally for small initial data in various functional spaces.
Before we address the types of initial data for which the existence of 
solutions has been established, we recall the scaling invariance of the equations.
If the pair $(u(x,t), p(x,t))$ solves the Navier-Stokes
equations \eqref{ns} in  ${\mathbb R}^{d}$ then $(u_{\lambda}(x,t), p_{\lambda}(x,t))$ with
$$ u_{\lambda}(x,t) = \lambda u(\lambda x, \lambda^2 t),$$
$$ p_{\lambda}(x,t) = \lambda^2 p(\lambda x, \lambda^2 t)$$
is a solution to the system \eqref{ns} for the initial data
\begin{equation} \label{scaluz} u_{0 \; \lambda} = \lambda u_0
  (\lambda x) \;\;. \end{equation}
The spaces which are
invariant under 
such a scaling are called critical spaces for the Navier-Stokes equations. Examples
of critical spaces for the Navier-Stokes in ${\mathbb R}^{d}$ are: 
\be{embed} 
 \dot{H}^{\frac{d}{2}-1} \hookrightarrow L^{d}  
\hookrightarrow  \dot{B}_{p|p < \infty, \infty}^{-1+\frac{d}{p}}  
\hookrightarrow BMO^{-1}.
\end{equation} 
The study of the Navier-Stokes equations in critical spaces was initiated by Kato \cite{K84}
and continued by many authors, see, for example, \cite{GM89}, \cite{T92}, \cite{C97}, 
\cite{P96}. In 2001 Koch and Tataru \cite{KT01} proved the existence of global solutions to 
\eqref{ns} - \eqref{0nsin} in  ${\mathbb R}^{d}$ corresponding to initial data small enough in 
$BMO^{-1}$. The space $BMO^{-1}$ has a special role since it is the largest 
critical space among the spaces listed in \eqref{embed} 
where such existence results are available. 

Motivated by the result of Koch and Tataru \cite{KT01} we analyze regularity
properties of the solution constructed by Koch and Tataru. More precisely, we 
show that under certain smallness condition of the initial data in $BMO^{-1}$, 
the solution $u$ to the Navier-Stokes equations \eqref{ns} - \eqref{0nsin} 
constructed in \cite{KT01} satisfies the following regularity property:
\be{GPSreg} 
t^{\frac{k}{2}} \nabla^k u \in  X^0, \mbox{ for all } k \in {\mathbb N} \cup \{0\}\, ,
\end{equation}
where $X^0$ denotes the space where the solution constructed by Koch and Tataru belongs 
(for a precise definition of $X^0$, see Section 2). 
In order to identify such a smoothing effect on solutions to the Navier-Stokes 
equations, we modify the argument \cite{KT01} in an appropriate way. 
The main tools in establishing the estimates sufficient to carry out 
the fixed point algorithm are properties of the Oseen kernel, modified maximal regularity estimates of the heat kernel, and a $TT^{*}$ argument. 
As a corollary we also prove that the  
solution to the Navier-Stokes equations constructed via such a fixed point
algorithm exhibits the decay in time of space derivatives of the type 
$$ \|\nabla^{k} u \|_{BMO^{-1}} \leq \frac{C}{t^{k/2}}, \mbox{ for } k \geq 1,$$
for any $t \geq 0$. 
 
Similar smoothing effects of solutions to the Navier-Stokes equations in 
the Lebesgue space $L^d$ were analyzed by Giga and Sawada \cite{GS02} 
and Dong and Du \cite{DD05}, and in the homogeneous Sobolev space $\dot{H}^{d/2 -1}$ 
by Sawada \cite{Sa03}. As in the case of this paper, the authors in \cite{GS02},
\cite{DD05} and \cite{Sa03} modified the corresponding existence results in order 
to identify more regular behavior of solutions to the Navier-Stokes equations 
and obtain consequently decay estimates. 

We note that long time behavior of solutions to the Navier-Stokes equations 
has been studied by applying a different approach too, based on splitting the initial
data in a certain critical space into a part that belongs to $L^2$ and
a part that is small in the corresponding critical space. We refer the reader, for example, to 
\cite{GIP02} and \cite{ADT04}, where the advantage of working with weak solutions have 
been extensively used, without identifying regularity estimates.
On the other hand, in order to study long time behavior, M. Schonbeck and 
her collaborators have successfully combined regularity estimates and decay 
estimates of solutions to the Navier-Stokes equations based on $L^2$ theory, 
see, for example, \cite{Sc95}. 

Often one obtains spatial analyticity of solutions to the Navier-Stokes equations as 
a consequence of the fixed point scheme used to establish existence and regularity of 
solutions to the Navier-Stokes equations. Here by spatial analyticity of solutions 
we mean convergence of the spatial Taylor series associated with the solution. 
Analysis of spatial analyticity of the Navier-Stokes equations was initiated  by Masuda 
\cite{M67}, and Kahane~\cite{Kah} while time analyticity was discussed by Foias and Temam \cite{FT89}. 
The study of analyticity of the Navier-Stokes equations was then continued by 
many authors. For example, Le Jan and 
Sznitman \cite{LJS97} considered solutions in a certain Besov space based on pseudomeasures, 
Gruji\'{c} and Kukavica \cite{GK98}, Lemari\'{e}-Rieusset \cite{L00}, 
Giga and Sawada \cite{GS02} studied spatial analyticity of the solutions to the 
Navier-Stokes equations in the critical Lebesgue space $L^d$, while Foias and Temam \cite{FT89} 
and Sawada \cite{Sa03} considered spatial analyticity for solutions in Sobolev spaces.  

With just a little bit of extra work, as a side-product of our main theorem we can obtain 
the spatial analyticity of the solution to the Navier-Stokes equations.

Finally, we apply our main regularity result for the Koch and Tataru solutions to 
(forward) self-similar solutions. We shall come back to this question in 
section~\ref{sectionselfsim}, but let us mention quickly what we obtain. Self-similar solutions
are given by non-linear profiles, and we get new regularity and decay estimates for these profiles.

To summarize in this paper we present regularity, decay and analyticity estimates for
solutions to the Navier-Stokes equations evolving from small initial data 
in $BMO^{-1}$. We do that by modifying the proof of existence of such solutions 
\cite{KT01} and without relying on $L^2$ theory. Also we present an application 
of our regularity result to self-similar solutions.

After the completion of the present article, we learned about the recent 
paper of Miura and Sawada \cite{M06} on the
regularity of Koch and Tataru solutions to the Navier-Stokes equations. 
However we emphasize the difference in the regularity result obtained
in \cite{M06} with respect to the regularity result that we obtain in this paper. 
More precisely, in \cite{M06} Miura and Sawada prove 
that the global Koch-Tataru solution to the system 
\eqref{ns} - \eqref{0nsin} evolving from small initial data in 
$BMO^{-1}$ (or a local in time solution in $VMO^{-1}$) 
satisfies the following regularity property:
\be{MSreg}
t^{\frac{k}{2}} \nabla^k u \in  N^0_{\infty}, 
\mbox{ for all } k \in {\mathbb N} \cup \{0\}\, .
\end{equation}
We remark that the spaces $X^0$ appearing in \eqref{GPSreg} and $N^0_{\infty}$  
are related via 
$$ \|u\|_{X^0} = \|u\|_{N^0_{\infty}} + \|u\|_{N^0_C},$$ 
(see Section 2 for the precise definition of spaces). Thus our regularity result 
can be understood as an extension of the result of Miura and Sawada. 
Indeed, a major part of our paper concentrates on obtaining the regularity result
for the Carleson part of the norm\footnote{We use this to obtain a regularity result for self-similar solutions too.}
given by $N^0_C$.
We remark that both papers give the analyticity result, obtained by slightly 
different means. Although we use the same pointwise bound on the heat kernel 
as Miura and Sawada, we do not use the Gronwall type inequality 
which has been obtained in earlier 
works \cite{GG} and \cite{GS02} (see also Lemma 2.7 in \cite{M06} ).

\subsection*{Organization of the paper} 
In Section 2 we introduce the notation that we shall use throughout the paper 
and state our results. In Section 3 we formulate three results in harmonic analysis
that we shall use in the proof of our main regularity theorem. Section 4 offers a proof of 
the main regularity result: Theorem \ref{reg}, as well as some consequences. 
In Section 5 we show the spatial analyticity result: Theorem \ref{analyt}. 
In Section 6 we prove a regularity result for self-similar solutions. 

\subsection*{Acknowledgments} We would like to express our thanks to 
Vladimir \v{S}ver\'{a}k for suggesting the problem, and to Hongjie Dong, Isabelle 
Gallagher, Patrick Gerard and Dong Li for interesting and useful discussions during 
the writing of this paper.

\section{Preliminaries and the statement of the result} 

\subsection{Notations and definition of the functional spaces} 

Throughout the paper we use $A \lesssim B$ to denote an estimate of the form $A\leq CB$ for some
absolute constant $C$. If $A \lesssim B$ and $B \lesssim A$ we write $A \sim B$. 

\bigskip

First we recall the definition of the space $BMO^{-1}$ as presented in the paper of 
Koch and Tataru \cite{KT01}:
\be{defbmo-1KT}
\|f(\cdot)\|_{BMO^{-1}} = \sup_{x_0,R} 
\left( \; \frac{1}{|B(x_0,\sqrt{R})|} \int_{0}^{R} \int_{B(x_0,\sqrt{R})}  
|e^{t \Delta} f(y)|^{2} \; dy \; dt \; \right)^{\frac{1}{2}}.
\end{equation} 

In  \cite{KT01} Koch and Tataru proved the following existence theorem  
for the solutions to the Navier-Stokes equations:
\begin{thm} \label{KTtheorem} 
The Navier-Stokes equations \eqref{ns} - \eqref{0nsin} with $f=0$ 
have a unique global solution in $X^0$ 
\begin{align*}
\| u(\cdot, \cdot) \|_{ X^0 } & = \sup_{t} t^{\frac{1}{2}} \|u(\cdot, t)\|_{L^{\infty}} \\
& + \sup_{x_0,R} \left( \; \frac{1}{|B(x_0,\sqrt{R})|} \int_{0}^{R} \int_{B(x_0,\sqrt{R})}  
|u(y,t)|^{2} \; dy \; dt \; \right)^{\frac{1}{2}},
\end{align*}
for all initial data $u_0$ with $\nabla \cdot u_0 = 0$ which are small in $BMO^{-1}$.
\end{thm}
We shall call such a solution the Koch-Tataru solution to the Navier-Stokes equations. 
It was pointed out by Auscher et al \cite{ADT04} that the Koch-Tataru solution actually verifies 
\begin{equation}
\label{ausch}
u \in L^{\infty}([0, \infty); BMO^{-1}) \,\,.
\end{equation}

Inspired by studying smoothing properties of the Koch-Tataru solutions of the
Navier-Stokes equations, for a nonnegative integer $k$ we introduce the space 
$X^k$ which is equipped with the following norm: 
$$ \|u\|_{X^k} = \|u\|_{N^k_{\infty}} + \|u\|_{N^k_C},$$
where
\begin{align*} 
& \|u\|_{N^k_{\infty}}  = \sup_{\alpha_1 + \dots + \alpha_d = k} \, \sup_{t} t^{\frac{k+1}{2}} 
\| \partial_{x_1}^{\alpha_1} \dots \partial_{x_d}^{\alpha_d}
 u(\cdot, t)\|_{L^{\infty}}, \\
& \|u\|_{N^k_C}  = \\
& \;\;\;\;\;\sup_{\alpha_1 + \dots + \alpha_d = k} \, \sup_{x_0,R} 
\left( \; \frac{1}{|B(x_0,\sqrt{R})|} \int_{0}^{R} \int_{B(x_0,\sqrt{R})} 
|t^{\frac{k}{2}} \partial_{x_1}^{\alpha_1} \dots \partial_{x_d}^{\alpha_d} 
u(y,t)|^{2} \; dy \; dt \; \right)^{\frac{1}{2}}\,\,.
\end{align*}

\medskip

{\bf{Remark}} In the following, we will denote  
$\nabla^{k} u = \partial_{x_1}^{\alpha_1} \dots \partial_{x_d}^{\alpha_d} u$ 
with $ (\alpha_1, \alpha_2, ..., \alpha_d) \in {\mathbb{N}}_{0}^{d}$ and 
$k = \alpha_1 + ... + \alpha_d$. More generally, if 
$\| \cdot \|$ is a norm, we will always write $\| \nabla^k f \|$ instead of 
$\displaystyle \sup_{\alpha_1 + \dots + \alpha_d = k} 
\| \partial_{x_1}^{\alpha_1} \dots \partial_{x_d}^{\alpha_d} f \|$.

\medskip

Hence
\be{property} 
u \in X^k \mbox{ if and only if } t^{\frac{k}{2}} \nabla^{k} u \in X^{0}. 
\end{equation} 

\subsection{Formulation of the results} 

Now we are ready to formulate the main result of this note: 

\begin{thm} \label{reg} 
There exists $\epsilon = \epsilon(d)$ such that if $\| u_0 \|_{BMO^{-1}} < \epsilon$, 
the Koch-Tataru solution $u$ associated to the initial value
problem \eqref{ns} - \eqref{0nsin} with $f=0$ verifies
$$
t^{\frac{k}{2}} \nabla^{k} u \in X^{0}
$$
for any $k \geq 0$.
\end{thm}


Theorem \ref{reg} implies the following decay in time of the space derivatives:

\begin{cor}
\label{tspderiv}
If $\|u_0\|_{BMO^{-1}} < \epsilon(d)$, the Koch-Tataru solution $u$ satisfies
\begin{equation}
\label{BMO-1dec}
\| \nabla^{k} u \|_{BMO^{-1}} \leq \frac{C}{t^{k/2}},
\end{equation}
for any $t\geq 0$ and any $k \geq 0$. 
\end{cor}

\medskip

{\bf{Remarks}}
\begin{itemize}
\item The case $k = 0$ in Corollary \ref{tspderiv} 
corresponds to the result \eqref{ausch} 
obtained by Auscher, Dubois and Tchamitchian in \cite{ADT04}. 

\item Auscher, Dubois and Tchamitchian proved in~\cite{ADT04} that any
  global solution with an initial condition\footnote{not necessarily 
small in $BMO^{-1}$} in the closure of the Schwartz class in $BMO^{-1}$ goes to $0$ in the
  $BMO^{-1}$ norm as $t$ goes to infinity. Combined with Corollary
  \ref{tspderiv}, this yields that any such solution satisfies the
  bounds~(\ref{BMO-1dec}) for $t$ large. 
Furthermore, the arguments given in the present article can easily be adapted
to handle the case of solutions defined on a finite time interval. 

As a consequence, we get that the bounds~(\ref{BMO-1dec}) are verified for any $t$
as soon as we have a global solution with an initial condition in the
  closure of the Schwartz class in $BMO^{-1}$.
\end{itemize}

Also the proof of Theorem \ref{reg} implies the following result: 

\begin{thm} \label{analyt} 
If $\|u_0\|_{BMO^{-1}} < \epsilon(d)$, then the Koch-Tataru global solution $u$ is 
space analytic. 
\end{thm}

\subsection{Application to self-similar solutions}

\label{sectionselfsim}

We recall that a solution to the Navier-Stokes equations \eqref{ns} - \eqref{nsdiv} with
$f=0$ is called self-similar if it can be written as  
\be{self}
u(x,t) = \lambda(t) \; U\left(\lambda(t)x\right), \; \; \; p(x,t) = \lambda^{2}(t) \; P\left(\lambda(t)x\right),
\end{equation} 
where $U(x)$ is a divergence-free vector field and $P(x)$ is a scalar field. In particular we distinguish two types
of self-similar solutions: 

\begin{enumerate} 
\item A backward self-similar solution of \eqref{ns} - \eqref{nsdiv} is a solution of the type \eqref{self} 
with $\lambda(t) = \frac{1}{\sqrt{2a(T-t)}}$, $a>0$, $T>0$ and $t < T$. In this case
the pair $(U,P)$ is a solution to the system:
\begin{align*}
& -\Delta U + a U + a(x \cdot \nabla)U + (U \cdot \nabla)U + \nabla P = 0 \\
& \nabla \cdot U = 0.
\end{align*} 

\item A forward self-similar solution of \eqref{ns}- \eqref{nsdiv} is a solution of the type \eqref{self} 
with $\lambda(t) = \frac{1}{\sqrt{2a(T+t)}}$, $a>0$, $T>0$ and $t > - T$. In this case
the pair $(U,P)$ is a solution to the system:
\begin{align}
& -\Delta U - a U - a(x \cdot \nabla)U + (U \cdot \nabla)U + \nabla P = 0 \label{selfF} \\
& \nabla \cdot U = 0. \nonumber
\end{align} 
\end{enumerate} 

It has been established that backward self-similar solutions sufficiently decaying 
at infinity do not exist, see, the work of Ne\v{c}as, Ru\v{z}i\v{c}ka and \v{S}ver\'{a}k 
\cite{NRS} and the work of Tsai \cite{Ts}. 

Various forward self-similar solutions were constructed 
provided the data is small in some critical space, see, for example, 
works of Giga and Miyakawa \cite{GM89}, Cannone, Meyer and Planchon \cite{C94} and
Cannone and Planchon \cite{C96}. Whenever non-trivial forward self-similar solutions exist, 
their initial data $u(x,0)$ is a homogeneous function of degree $-1$. Hence 
it is natural to analyze self-similar solutions in critical spaces containing 
homogeneous functions of degree $-1$. 

Now let us suppose that $u_0$ is self similar, i.e. 
\be{iselfsim} 
u_0(x) = \lambda u_0 (\lambda x) \;\;\;\; \mbox{for any}\;\; \lambda > 0 \,\,.
\end{equation} 
Then if some uniqueness property is verified $u$ itself will be self-similar, 
in the sense that
\begin{equation}
\label{selfsim}
u(x,t) = \lambda u ( \lambda x , \lambda^2 t ) \;\;\;\;\;\;\mbox{for any}\;\; \lambda > 0 \,\,.
\end{equation}
At least formally, this gives the existence of a profile $\psi$ such that
$$
u(x,t) = \frac{1}{\sqrt{t}} \psi \left( \frac{x}{\sqrt{t}} \right) \,\,
$$
which satisfies the elliptic equation \eqref{selfF} with  $a = 1/2$ and $T=0$ i.e.
$$
- \Delta \psi - \frac{1}{2} \psi 
- \frac{1}{2} (y \cdot \nabla) \psi
+ (\psi \cdot \nabla) \psi 
+ \nabla P = 0
$$
for some scalar function $P$.

Recently, Gruji\'{c} \cite{G06} proved regularity of a self-similar 
solution in $L^{\infty}_t L^{3}_x$ associated to small data in 
the Lebesgue space $L^3_x$. 
With the help of Theorem \ref{reg}, we are able to generalize this result 
to small data in $BMO^{-1}$. More precisely, we prove: 

\begin{thm}
\label{thselfsim}
There exists $\epsilon, \delta >0$ such that:
if $u_0$ is self similar (equation~(\ref{iselfsim}) is verified) with a norm in $BMO^{-1}$ 
smaller than $\epsilon$, then there exists a unique solution $u$ to the Navier-Stokes equations 
\eqref{ns} - \eqref{nsdiv} with a norm in $X^0$ smaller than $\delta$. This solution is 
given by a profile $\psi$:
$$
u(x,t) = \frac{1}{\sqrt{t}} \psi \left( \frac{x}{\sqrt{t}} \right) \,\,.
$$

Also for any nonnegative integer $k$ we have
\begin{equation}
\label{estimselfsim}
\nabla^k \psi \in L^\infty \;\;\;\;\mbox{and}\;\;\;\;
\int_{{\mathbb{R}}^d} \frac{1}{|y|^d} \left| \nabla^k \psi(y) \right|^2 \,dy < \infty \,\,.
\end{equation}
\end{thm}

The first part of~(\ref{estimselfsim}) implies that $\psi \in \mathcal{C}^\infty$; 
the second part gives in particular that all derivatives of $\psi$
vanish at $0$, and an indication on their decay at infinity.

This theorem is proved in section~\ref{proofselfsim}.

\section{Three results of harmonic analysis}

\subsection{A Carleson-type estimate}

We shall present a modification of the result given in Lemari\'{e}'s book \cite{L02}
Lemma 16.2, page 163. The result of Lemma 16.2 from \cite{L02} was originally proved by 
Koch and Tataru \cite{KT01}. A modified version will be used in Section \ref{pro} 
in order to obtain a certain bound on the nonlinear term. 

\begin{lem}\label{mLe16.2}
For $N(x,t)$ defined on ${\mathbb{R}}^{d} \times (0,1)$, let $A(N)$ be the
quantity
$$A(N) =\sup_{x_0 \in {\mathbb{R}}^{d}} \sup_{0 < t < 1} t^{-\frac{d}{2}}
\int_{0}^{t} \int_{|x-x_0| < \sqrt{t}} |N(x,s)| \; dx \; ds.$$
Then there exists a constant $b(k)$ such that the following inequality holds for 
$\beta_{k}(x,t) = 
t^{\frac{k}{2}} (-\Delta)^{\frac{k+1}{2}} e^{t\Delta} \int_{0}^{t} N(x,s) \; ds$:
$$ \int_{0}^{1} \int_{{\mathbb{R}}^{d}} |\beta_{k} (x,t)|^2 \; dx \; dt \leq 
b(k) A(N) \int_{0}^{1} \int_{{\mathbb{R}}^{d}} |N(x,s)| \; dx \; ds.$$
\end{lem}

\begin{proof} 
\begin{align} 
& \int_{0}^{1} \int_{{\mathbb{R}}^{d}} |\beta_k (x,t)|^2 \; dx \; dt \nonumber \\
& = \int_{0}^{1} \langle
\int_{0}^{t} t^{\frac{k}{2}} (-\Delta)^{\frac{k+1}{2}} e^{t\Delta} N(\cdot\;, s) ds, 
\int_{0}^{t} t^{\frac{k}{2}} (-\Delta)^{\frac{k+1}{2}} e^{t\Delta} N(\cdot\;, \sigma) d\sigma 
\rangle_{L^2(dx)} \; dt \nonumber \\
& = 2 Re \int \int \int_{0 < \sigma < s < t < 1} 
\langle
 N(\cdot\;, s), 
 t^{k} (-\Delta)^{k+1} e^{2t\Delta} N(\cdot\;, \sigma) 
\rangle_{L^2(dx)} \; d\sigma \; ds \; dt \nonumber \\
& = 2 Re \int \int_{0 < \sigma < s < 1} 
\langle N(\cdot\;, s), 
\left( \int_{s}^{1} t^k (-\Delta)^{k+1} e^{2t\Delta} \; dt \right)
\; N(\cdot\;, \sigma) 
\rangle_{L^2(dx)} \; d\sigma \; ds \nonumber \\
& = 2 Re \int_{0}^{1}   
\langle  N(\cdot\;, s), \int_{0}^{s} \left( L_k(1) - L_k(s) \right) 
N(\cdot\;, \sigma) \; d\sigma \rangle_{L^2(dx)}\; ds, \label{intL} 
\end{align} 
where $L_k$ is given via the following integral 
$$ \int t^k (-\Delta)^{k+1} e^{2t\Delta} \; dt = L_k(t) + C,$$
which after performing a sequence of integration by parts gives the following operator  
\be{calcL} 
L_k(t) = \sum_{m=0}^{k} b_m(k) t^m (-\Delta)^{m} e^{2t \Delta}, 
\end{equation} 
with the constant $b_m(k)$ depending on both $m$ and $k$; however such a dependence will not 
be crucial for the argument, so it is not specified here. 

Let 
$$\gamma_k = \sup_{0 < s \leq 1} 
\| \int_{0}^{s} L_k(s) N(\cdot\;, \sigma) \; d\sigma \|_{L^{\infty}_{dx}}.$$
Then \eqref{intL} implies that 
\be{bebound} 
 \int_{0}^{1} \int_{{\mathbb{R}}^{d}} |\beta_k (x,t)|^2 \; dx \; dt 
\lesssim \gamma_k \int_{0}^{1} \int_{{\mathbb{R}}^d} |N(x, s)|\; dx \; ds.
\end{equation} 

In order to conclude the proof we shall obtain an upper bound on $\gamma_k$. 
We follow the approach of Koch and Tataru \cite{KT01} as presented in Lemarie's book
\cite{L02}. More precisely, let us write the kernel of the operator $L_k$ as
$ \frac{1}{s^{d/2}} W_k(\frac{x}{\sqrt{s}})$ with a positive $W_k \in {\mathcal S}$. 
Then we have: 
\begin{align*} 
& \left|\int_{0}^{s} L_k(s)N(\cdot\;, \sigma) \; d\sigma \right| \\
& \leq b(k) \int_{0}^{s} \int_{ {\mathbb R}^d } \frac{1}{s^{d/2}} W_k(\frac{x-y}{\sqrt{s}}) 
|N(y, \sigma)| \; dy \; d\sigma \\
& = \sum_{q \in {\mathbb Z}^d} \int_{0}^{s} \int_{x - \sqrt{s} (q + [0,1]^d)} 
\frac{1}{s^{d/2}} W_k(\frac{x-y}{\sqrt{s}}) |N(y, \sigma)| \; dy \; d\sigma \\ 
& \leq \sum_{q \in {\mathbb Z}^d} \sup_{z \in q + [0,1]^d} W_k(z) \; 
\frac{1}{s^{d/2}} \int_{0}^{s} \int_{x - \sqrt{s}(q+[0,1]^d)}
|N(y, \sigma)| \; dy \; d\sigma, 
\end{align*} 
which implies that 
\be{gabound} 
\gamma_k \lesssim b(k) A(N). 
\end{equation} 
Now we combine \eqref{bebound} and \eqref{gabound} 
to conclude that 
$$  \int_{0}^{1} \int_{{\mathbb{R}}^{d}} |\beta_k (x,t)|^2 \; dx \; dt 
\lesssim b(k) A(N) \int_{0}^{1} \int_{{\mathbb{R}}^d} |N(x, s)|\; dx \; ds,$$
and the lemma is proved. 
\end{proof}

\subsection{Generalized maximal regularity of the heat kernel}

We will need in the following section 
a generalization of the maximal regularity of the heat kernel.

\begin{prop}
\label{laproposition}
If $r$ is a natural number, the operators
$$
P_r : f \mapsto \int_0^t e^{(t-s)\Delta} (t-s)^r \Delta^{r+1} f(s) \; ds
$$
and
$$
Q_r : f \mapsto \int_0^t e^{(t-s)\Delta} (t-s)^r (\sqrt{t} - \sqrt{s}) \Delta^{r+1} \sqrt{-\Delta} f(s) \; ds
$$
are bounded on $L^2([0,T], L^2({\mathbb{R}}^d))$ for any $T \in [0,\infty]$ with constants respectively 
$p(r)$ and $q(r)$.
\end{prop}

{\bf{Remark}} The classical maximal regularity of the heat kernel corresponds 
to the boundedness of $P_0$ on $L^p ([0,T],L^q({\mathbb{R}}^d))$.

\begin{proof} We prove the theorem for $T=\infty$ ; the result for any other $T$ is then an easy consequence.
 We also omit to write the dependence of the constants on $r$.

Let us begin with the boundedness of $P_r$, which can be dealt with 
simply by using the Fourier transform. Indeed, taking the Fourier transform in 
space ${\mathcal{F}}_x$, $P_r$ becomes
$$
{\mathcal{F}}_x P_r f = \int_0^t e^{-(t-s)|\xi|^2} (t-s)^r |\xi|^{2r+2} 
{\mathcal{F}}_x f(s) \; ds \,\,.
$$
If we assume that $f$ is zero for $t<0$, the above expression is nothing but 
the convolution of ${\mathcal{F}}_x f$ by $\chi(t) e^{-t|\xi|^2} t^r |\xi|^{2r+2}$, 
where $\chi$ is the characteristic function of $(0, \infty)$. Denoting by
${\mathcal{F}}_t$ the time Fourier transform, we see that $P_r$ is a space-time 
Fourier multiplier of symbol ${\mathcal{F}}_t \left[ \chi(t) e^{-t|\xi|^2} t^r |\xi|^{2r+2} \right]$. 
Therefore, $P_r$ will be bounded on $L^2$ if and only if
$$
{\mathcal{F}}_t \left[ \chi(t) e^{-t|\xi|^2} t^r |\xi|^{2r+2} \right]
$$
is bounded. Removing the dilation by $|\xi|^2$ in the above expression, we see that 
it suffices to show that
$$
{\mathcal{F}}_t \left[ \chi(t) e^{-t} t^r \right]
$$
is a bounded function. But this is the case since $\chi(t) e^{-t} t^r$ belongs to $L^1$. 
This proves the boundedness of $P_r$ on $L^2$.

\bigskip

The case of $Q_r$ is more delicate because, due to the $(\sqrt{t} - \sqrt{s})$ factor, $Q_r$ is not a convolution operator, and therefore we cannot use the Fourier transform as easily as above.

Using the space Fourier transform, $Q_r f$ can be expressed as
\begin{equation*}
\begin{split}
{\mathcal{F}}_x Q_r f(t)  & =  \int_0^t (t-s)^r (\sqrt{t} - \sqrt{s}) 
e^{(s-t)|\xi|^2} |\xi|^{2r+3} 
{\mathcal{F}}_x f(s,\xi) \; ds \\
& \overset{\mbox{def}}{=} \int_0^t \alpha(\xi,s,t) {\mathcal{F}}_x f(s,\xi) \; ds
\end{split}
\end{equation*}

To conclude, we use the following lemma

\begin{lem}
\label{petitlemme}
The operator
$$
A : g \mapsto \int_0^\infty \alpha(1, s,t) g(s) ds
$$
is bounded on $L^2(\mathbb{R})$.
\end{lem}

Before proving this lemma, let us show how it will enable us to complete the proof of Proposition~\ref{laproposition}. By homogeneity, Lemma~\ref{petitlemme} shows that the operator
\begin{equation*}
\begin{split}
A_{\xi} : & L^2({\mathbb{R}}) \longrightarrow  L^2(\mathbb{R})\\
& g  \mapsto  \int_0^\infty \alpha(\xi,s,t) g(s) ds
\end{split}
\end{equation*}
is bounded with a bound independent of $\xi$. So we get
\begin{equation*}
\begin{split}
\|Q_r f\|_{L^2({\mathbb{R}}_+, L^2({\mathbb{R}}^d))} 
& = C \| {\mathcal{F}}_x Q_r f\|_{L^2({\mathbb{R}}_+, L^2({\mathbb{R}}^d))} \\
& = C \| \| A_{\xi} {\mathcal{F}}_x f(s,\xi) \|_{L^2(\mathbb{R}_+)} \|_{L^2({\mathbb{R}}^d)} \\
& \leq C \| \| {\mathcal{F}}_x f(s,\xi) \|_{L^2(\mathbb{R}_+)} \|_{L^2({\mathbb{R}}^d)}\\
& = C \| f\|_{ L^2( {\mathbb{R}}_{+}, L^2({\mathbb{R}}^d) ) } \,\,,
\end{split}
\end{equation*}
which proves the proposition. \end{proof} 

{\bf{Proof of Lemma~\ref{petitlemme}}}
We want to prove that the operator
$$
A : f \mapsto \int_0^t \alpha(s,t) f(s) \,ds
$$
with
$$
\alpha(s,t) = (t-s)^r (\sqrt{t} - \sqrt{s}) e^{(s-t)}
$$
is bounded on $L^2(\mathbb{R}_+)$. The kernel $\alpha$ is non-negative, and, 
due to the inequality $\sqrt{t}-\sqrt{s} \leq \sqrt{t-s}$ for $0<s<t$, we have the majorization
$$
\alpha(s,t) \leq (t-s)^r \sqrt{t-s} e^{(s-t)} \,\,.
$$
Using the Fourier transform in time, it is easy to see that this last kernel defines a bounded operator 
on $L^2(\mathbb{R}_+)$. This implies that $A$ is also a bounded operator on $L^2(\mathbb{R}_+)$.

\subsection{The Oseen kernel} In the proof of the main theorem we will use
the following proposition about the Oseen kernel whose proof can be found 
in \cite{L02}, Proposition 11.1.

\begin{prop} \label{Oseen} 
The integral kernel of the operator $\nabla^{k+1}{\mathbb{P}} e^{t\Delta}$ (where $\mathbb{P}$ represents the Leray projection onto the divergence-free vector fields)
is bounded pointwise by
$$
\frac{K(k)}{\left(\sqrt{t}+|x|\right)^{d+k+1}} \,\,.
$$
\end{prop}

Although the above proposition is enough in order to obtain the
regularity result stated in Theorem \ref{reg}, we recall the proposition
given by Miura and Sawada \cite{M06} which expresses a slightly different 
pointwise bound on the derivatve of the heat kernel (note 
the difference in the decay).

\begin{prop} \label{OseenMS} 
The integral kernel of the operator $\nabla^{k+1}{\mathbb{P}} e^{t\Delta}$ (where $\mathbb{P}$ represents the Leray projection onto the divergence-free vector fields)
is bounded pointwise by
$$
\frac{C^k k^{k/2}}{t^{k/2}\left(\sqrt{\frac{t}{k}} + |x|\right)^{d+1}} \,\,.
$$
\end{prop}

Also let us recall here the following boundedness property of the
heat kernel (see, for example, \cite{KF}): 
\begin{equation} \label{heatb} 
\|\nabla e^{t \Delta} u \|_{L^{\infty}} \leq \frac{C}{\sqrt{t}} \|u\|_{L^{\infty}}. 
\end{equation}

\setcounter{equation}{0}
\section{Proof of the main theorem} \label{pro}

In \cite{KT01} Koch and Tataru proved that if the initial data $u(x,0)$ are divergence-free and 
sufficiently small in $BMO^{-1}$, then there exists a solution $u(x,t) \in X^{0}$ to 
the integral Navier-Stokes equations 
\be{intNS} 
u(x,t) = e^{t \Delta} u(x,0) - B(u,u)(x,t),
\end{equation} 
where 
\be{B} 
B(u,v)(x,t) = \int_{0}^{t} e^{(t-s)\Delta} {\mathbb{P}} \nabla \cdot (u(x,s) \times v(x,s)) \; ds,
\end{equation}
where $u \times v$ denotes the tensor product of $u$ and $v$.

Let us fix a positive integer $k$. We shall prove that $u(x,t)$ given by \eqref{intNS} 
is in  $X^{k}$ if the data is small enough.

\subsection{Linear term} 

\begin{prop}
\label{linear}
For any $k \geq 0$, there exists a constant $C(k)$ such that
$$
\| e^{t \Delta} u_0 \|_{X^k} \leq C(k) \|u_0\|_{BMO^{-1}} \,\,.
$$
\end{prop}

\begin{proof}
Let us begin with the $L^\infty$ part of the norm. We would like to show that
\begin{equation}
\label{decinf}
\| e^{t \Delta} \nabla^{k} u_0 \|_{\infty} \leq C(k) \|u_0\|_{BMO^{-1}} t^{-\frac{k+1}{2}} \,\,.
\end{equation}
This is a consequence of the three following facts, whose proofs can be found in the book of Lemari\'e~\cite{L02} (the same reference also provides a 
definition of the Besov spaces $\dot{B}_{\infty,\infty}^{-l}$).
\begin{itemize}
\item The estimate 
$\displaystyle \| e^{t \Delta} f \|_{\infty} \leq C(k) t^{-\frac{k+1}{2}}$ 
holds if and only if $f \in 
\dot{B}_{\infty,\infty}^{-k-1}$.
\item The space $BMO^{-1}$ embeds continuously into $\dot{B}_{\infty,\infty}^{-1}$.
\item The operator $\nabla$ is bounded from $\dot{B}_{\infty,\infty}^{-l}$ to $\dot{B}_{\infty,\infty}^{-l-1}$, for any $l \geq 0$.
\end{itemize}

Let us now turn to the Carleson part of the norm. We would like to prove that
$$
\sup_{x_0,R} \; \frac{1}{|B(x_0,\sqrt{R})|} 
\int_{0}^{R} \int_{B(x_0,\sqrt{R})}  
|t^{\frac{k}{2}} \nabla^{k} e^{t\Delta} u_0(y) |^{2} \; dy \; dt \; \leq C(k) \|u_0\|_{BMO^{-1}}^2 \,\,.
$$
Since $u_0 \in BMO^{-1}$, it can be written as a finite sum of derivatives of functions 
in $BMO$, $u_0 = \sum_i \partial_{x_i} f_i$, with 
$\sum \| f_i \|_{BMO} \sim C \|u_0\|_{BMO^{-1}}$. We assume for simplicity that $u_0 = \partial_{x_1} f$. The above inequality becomes
$$
\sup_{x_0,R} \; \frac{1}{|B(x_0,\sqrt{R})|} 
\int_{0}^{R} \int_{B(x_0,\sqrt{R})}  
|t^{\frac{k}{2}} \nabla^{k} e^{t\Delta} \partial_{x_1} f(y) |^{2} \; dy \; dt \; \leq \|f\|_{BMO}^2 \,\,.
$$
Let us now denote by ${\psi}^k$ the inverse Fourier transform of 
$(i\xi)^k i \xi_1 e^{-|\xi|^2}$, and 
$\displaystyle \psi_t^k(\cdot) 
= \frac{1}{t^{d/2}} {\psi}^k \left( \frac{\cdot}{\sqrt{t}} \right)$. 
Performing additionnally the change of variables 
$s = \sqrt{t}$, the above inequality can be rewritten as
$$
\sup_{x_0,{R}} \; \frac{1}{|B(x_0,\sqrt{R})|} 
\int_{0}^{\sqrt{R}} \int_{B(x_0,\sqrt{R})}  
| \psi_s^k * f |^{2} \; dy \; \frac{ds}{s} \; \leq C(k) \|f\|_{BMO}^2 \,\,.
$$
This last inequality holds true: this is one of the possible definitions of $BMO$, see Stein~\cite{S93}, Chapter 4.
\end{proof}

\subsection{Nonlinear term: the main estimate}

In order to simplify the notation let us for $k \geq 0$ denote by $\widetilde{X}^k$ the space
$$ \widetilde{X}^k = \cap_{l=0}^{k} X^l.$$
equipped with the norm $\sum_{l=0}^{k} \|\cdot\|_{X^l}$.

We shall prove that the bilinear operator maps:
\be{wantbil} 
B: \widetilde{X}^k \times \widetilde{X}^k \rightarrow \widetilde{X}^k \,\,.
\end{equation} 

More precisely, we have the following proposition

\begin{prop}
\label{estimb}
For any $k \geq 1$,
\begin{align} 
\begin{split} \label{particular}
\|B(u,v)\|_{X^k} 
& \leq C_0(k) \; \|u\|_{X^0} \|v\|_{X^0} \\ 
& + C_1 \; \|u\|_{X^0} \|v\|_{X^k} + C_1 \; \|v\|_{X^0} \|u\|_{X^k} \\
& + C(k) \; \|u\|_{\widetilde{X}^{k-1}} \|v\|_{\widetilde{X}^{k-1}} \,\,,
\end{split} 
\end{align} 
where the constant $C_1$ does not depend on $k$.
\end{prop} 

In the case when $k=0$ the following estimate has been proved
by Koch and Tataru in \cite{KT01}.
\begin{equation} 
\label{particular0} 
\|B(u,v)\|_{X^0} \leq C \|u\|_{X^0} \|v\|_{X^0}. 
\end{equation}

Subsections~\ref{nkinorm} and~\ref{ncknorm} are devoted to the proof of Proposition 
\ref{estimb}. In particular, details of the proof of \eqref{particular} are given. 
By using the convention that $\|u\|_{\widetilde{X}^{-1}} = \|v\|_{\widetilde{X}^{-1}} = 0$
in the proof of \eqref{particular}, a proof of \eqref{particular0} 
follows too. Thus in Subsections~\ref{nkinorm} and~\ref{ncknorm} we do not 
distinguish the case $k=0$.

\subsection{Nonlinear term : the $N^k_\infty$ norm}

\label{nkinorm}

Here, we shall prove that 
\begin{align} 
\begin{split} \label{BNinf}
\|B(u,v)\|_{N^k_{\infty}} 
& \leq c_0(k) \; \|u\|_{X^0} \|v\|_{X^0} \\ 
& + c_1 \; \|u\|_{X^0} \|v\|_{X^k} + c_1 \; \|v\|_{X^0} \|u\|_{X^k} \\
& + c(k) \; \sum_{l=1}^{k-1} \|u\|_{ N^{l}_{\infty} } \|v\|_{ N^{k-l}_{\infty} },
\end{split} 
\end{align} 
where the constant $c_1$ does not depend on $k$.

If $0 < s < t(1 - \frac{1}{m})$ (with $m=m(k)$ to be determined later)
 we use the estimate on the Oseen kernel Proposition \ref{OseenMS} 
to obtain
\begin{align*}
&  \int_{0}^{t(1 - \frac{1}{m})} 
|\nabla^{k} e^{(t-s)\Delta} {\mathbb P} \nabla \cdot (u(x,s) \times v(x,s)) | \; ds\\
& \leq C^k k^{\frac{k}{2}} \; \int_{0}^{t(1- \frac{1}{m})} 
\int_{{\mathbb R}^{d}} \frac {1} { (t-s)^{k/2} \left( \sqrt { \frac{t-s}{k}} + |x-y|\right)^{d+1}} |u (y,s)| |v(y,s)| \; dy \;ds\\  
& \lesssim  C^k k^{\frac{k}{2}} \left( \frac{m}{t}\right)^{(d+k+1)/2} \; \int_{0}^{t(1-\frac{1}{m})} \int_{{\mathbb R}^{d}} \frac {1}{\left( \frac{1}{\sqrt{k}} + 
\frac{|x-y|}{\sqrt{t-s}} \right)^{d+1}} |u (y,s)| |v(y,s)| \; dy \;ds\\ 
& \lesssim  C^k k^{\frac{k}{2}} \left( \frac{m}{t}\right)^{(d+k+1)/2} \; \int_{0}^{t(1-\frac{1}{m})} \sum_{q \in {\mathbb Z}^{d}} 
\frac{\int_{x-y \in \sqrt{t} (q+[0,1]^{d})} |u (y,s)| |v(y,s)| \; dy }
{ \left(\frac{1}{\sqrt{k}} + |q| \right)^{d+1} } \; ds, \label{infsum} 
\end{align*} 
which after applying the Cauchy-Schwartz inequality and using 
$$
\sum_{q \in {\mathbb Z}^{d}} \frac{1}{ \left(\frac{1}{\sqrt{k}} + |q| \right)^{d+1} } \sim \sqrt{k}
$$
implies
\begin{equation} \label{infint1} 
 \int_{0}^{t(1-\frac{1}{m})} |e^{(t-s)\Delta} {\mathbb P} \nabla^{k+1} \cdot (u(x,s) \times v(x,s)) | \; ds 
\leq  t^{ -\frac{k+1}{2} } c_0(k) \;  \|u\|_{X^0} \|v\|_{X^0},
\end{equation} 
with
\begin{equation} \label{ms-cons} 
c_0(k) =  C^k \; k^{\frac{k+1}{2}} \; m^{\frac{d+k+1}{2}}.
\end{equation}

If $t(1-\frac{1}{m}) \leq s < t$ we use \eqref{heatb} to obtain  
\begin{align*}
& | \nabla^k e^{(t-s)\Delta} {\mathbb P} \nabla \cdot (u(x,s) \times v(x,s)) | \\
& \lesssim \frac{1}{(t-s)^{1/2}} \sum_{l=0}^{k} 
\binom{k}{l}
\|\nabla^l u(\cdot, s)\|_{L^{\infty}} \|\nabla^{k-l} v(\cdot, s)\|_{L^{\infty}} \\ 
& \lesssim \frac{1}{(t-s)^{1/2}}  \sum_{l=0}^{k} 
\frac{1}{s^{\frac{l+1}{2} + \frac{k-l+1}{2}}}
\binom{k}{l} 
\|u\|_{ N^{l}_{\infty} } \|v\|_{ N^{k-l}_{\infty} }. 
\end{align*} 
Therefore 
\begin{align} 
& |\int_{t(1-\frac{1}{m})}^{t} \nabla e^{(t-s)\Delta} 
{\mathbb P} \nabla^{k+1} \cdot (u(x,s) \times v(x,s)) \; ds | \nonumber \\ 
& \lesssim \int_{t(1-\frac{1}{m})}^{t} \frac{1}{(t-s)^{1/2}} \frac{1}{s^{\frac{k+2}{2}}} \; ds \; 
 \sum_{l=0}^{k} \binom{k}{l} \|u\|_{ N^{l}_{\infty} } \|v\|_{ N^{k-l}_{\infty} } \nonumber \\
& \lesssim I(k,m,t) \left( \|u\|_{N^{0}_{\infty}} \|v\|_{N^{k}_{\infty}} 
+  \|v\|_{N^{0}_{\infty}} \|u\|_{N^{k}_{\infty}} + 
\sum_{l=1}^{k-1} \binom{k}{l} 
 \|u\|_{ N^{l}_{\infty} } \|v\|_{ N^{k-l}_{\infty} } \right), \label{infsplit} 
\end{align} 
where 
$$ I(k,m,t) =  \int_{t(1-\frac{1}{m})}^{t} \frac{1}{(t-s)^{1/2}} \frac{1}{s^{\frac{k+2}{2}}} \; ds.$$

After performing a change of variable $s = zt$, 
the integral $I(k,m,t)$ can be bounded from above as follows
\begin{align} 
I(k,m,t) 
& = t^{-\frac{k+1}{2}} \int_{1-\frac{1}{m}}^{1} 
\frac{1}{ (1-z)^{1/2} } \frac{1}{ z^{\frac{k+2}{2}} } dz \nonumber \\
& \leq t^{-\frac{k+1}{2}} \; (1 - \frac{1}{m})^{-\frac{k+2}{2}} 
\int_{1-\frac{1}{m}}^{1} (1-z)^{-1/2} dz \nonumber \\ 
& = 2 \;  t^{-\frac{k+1}{2}} g(m)\,\,,
\end{align}  
with 
$$g(m) = (1 - \frac{1}{m})^{-\frac{k+2}{2}} \frac{1}{\sqrt{m}}.$$ 
Now let\footnote{Such a choice of $m$
is motivated by the requirement $k^{\frac{k+1}{2}} m^{\frac{d+k+1}{2}} = k^{k-1}$, 
which in turn implies that the constant $c_0(k)$ appearing in 
\eqref{ms-cons} is like $C^{k-2} k^{k-1}$. This 
constant will be relevant in Section 5, where we prove the
analyticity of the solution.} 
\begin{equation} \label{them} 
m = m(k) = k^{\frac{k-3}{d+k+1}}.
\end{equation}

One can verify that $\lim_{k \rightarrow \infty} g\left(m(k)\right) = 0$, 
which implies that there exists $c$ such that for any integer $k$
$$|g\left(m(k)\right)| \leq c,$$
with a constant $c$ independent of $k$.  
Therefore  
\begin{equation} \label{Ibdd} 
I(k,m(k),t) \leq c t^{-\frac{k+1}{2}}, \; \; \mbox{ for all } k \geq 1.
\end{equation} 

Now the claim \eqref{BNinf} follows from \eqref{infint1}, \eqref{infsplit} and \eqref{Ibdd}.

\newpage

\subsection{Non-linear term: the $N_C^k$ norm}

\label{ncknorm}

\subsubsection{Splitting of $B$}

Now we shall prove that
\begin{align} 
\begin{split} \label{BNc}  
\|B(u,v)\|_{N^{k}_C}  
& \leq d_0(k) \; \|u\|_{X^0} \|v\|_{X^0} \\ 
& + d_1 \; \|u\|_{X^0} \|v\|_{X^k} + d_1 \; \|v\|_{X^0} \|u\|_{X^k} \\
& + d(k) \; \sum_{l=1}^{k-1} \|u\|_{ N^{l}_{\infty} } \|v\|_{ N^{k-l}_{\infty} },
\end{split} 
\end{align} 
where the constant $d_1$ does not depend on $k$. 


We split  $B(u,v)(x,t)$ as follows:
\be{Bsplit}
B(u,v)(x,t) = B_1 + B_2,
\end{equation} 
with
\begin{align*} 
B_1 & = \int_{0}^{t} e^{(t-s)\Delta} {\mathbb{P}} \nabla \cdot 
((1-\chi_{\sqrt{R},x_0}) u(x,s) \times v(x,s)) \; ds \\
B_2 & = \int_{0}^{t} e^{(t-s)\Delta} {\mathbb{P}} \nabla \cdot 
(\chi_{\sqrt{R},x_0} u(x,s) \times v(x,s)) \; ds,
\end{align*} 
where $\displaystyle \chi_{\sqrt{R},x_0} = \chi\left(\frac{x-x_0}{\sqrt{R}}\right)$ for a smooth function $\chi$ supported in $B(0,15)$ and equal to 1 on $B(0,10)$.

Also we further split $B_2$ as 
$$B_2 = B_{2}^{1} + B_{2}^{2},$$
where 
\begin{align*} 
B_{2}^{1} & = \frac{1}{\sqrt{-\Delta}} {\mathbb{P}} \nabla \cdot \int_{0}^{t} 
e^{(t-s)\Delta} \frac{\Delta}{\sqrt{-\Delta}} (Id-e^{s\Delta})
(\chi_{\sqrt{R},x_0} u(x,s) \times v(x,s)) \; ds, \\
B_{2}^{2} & = \frac{1}{\sqrt{-\Delta}} {\mathbb{P}} \nabla \cdot \sqrt{-\Delta} e^{t \Delta} 
\int_{0}^{t}  (\chi_{\sqrt{R},x_0} u(x,s) \times v(x,s)) \; ds.
\end{align*}

\subsubsection{Estimate for $B_1$}

To estimate $B_1$ we use a similar approach to the one that was applied to obtain the estimate 
\eqref{infint1}. However now we use the bound on the Oseen kernel expressed in Proposition \ref{Oseen}. 
More precisely, since $t < R$ we have 
\begin{align*}
& |t^{\frac{k}{2}} \nabla^k B_1(x,t)| \\
& \lesssim K(k) t^{\frac{k}{2}} 
\int_{0}^{t} \int_{|y-x_0| \geq 10 \sqrt{R}} 
\frac{1}{(\sqrt{t-s} + |x - y|)^{d+k+1}}  |u(y,s)| |v(y,s)| \; dy \; ds \\
& \lesssim K(k) R^{\frac{k}{2}}
\int_{0}^{R} \int_{|y-x| \geq 10 \sqrt{R}} 
\frac{1}{ R^{\frac{d+k+1}{2}} 
\left( \frac{\sqrt{t-s}}{\sqrt{R}} + \frac{|x-y|}{\sqrt{R}} \right)^{d+k+1} }
|u(y,s)| |v(y,s)| \; dy \; ds \\
& \lesssim K(k) R^{\frac{k}{2}}
\int_{0}^{R} \sum_{q \in {\mathbb Z}^d } 
\frac{\int_{x-y \in \sqrt{R} (q+[0,1]^{d})} |u (y,s)| |v(y,s)| \; dy }
{ R^{\frac{d+k+1}{2}} \; |q|^{d+k+1} } \; ds \\ 
& \lesssim K(k) D(k) 
R^{ \frac{k}{2} - \frac{d+k+1}{2} + \frac{d}{2} } \| u\|_{X^0} \|v\|_{X^0} \\
& =  K(k) D(k) R^{-\frac{1}{2}}  \|u\|_{X^0} \|v\|_{X^0},
\end{align*} 
where $D(k)$ denotes the constant coming from the summation in $q$ and 
$K(k)$ denotes the kernel bound.

Hence
\be{estB1}
\frac{1}{R^{\frac{d}{2}}} 
\int_{0}^{R} \int_{B(x_0,\sqrt{R})} |t^{\frac{k}{2}} \nabla^k B_1(u,v)(x,t)|^{2} \; dx \; dt
\lesssim \left( K(k)D(k) \right)^2 \|u\|_{X^0}^{2} \|v\|_{X^0}^{2}.
\end{equation}

\subsubsection{Estimate for $B_2^1$}

We would like to estimate the $L^2$ norm on a parabolic cylinder of
$$
t^{k/2} \nabla^k B_2^1 = t^{k/2} \nabla^k \frac{{\mathbb{P}} \nabla}{\sqrt{-\Delta}} \cdot \int_{0}^{t} 
e^{(t-s)\Delta} \Delta \frac{Id-e^{s\Delta}}{\sqrt{-\Delta}} 
(\chi_{\sqrt{R},x_0} u(s) \times v(s)) \; ds \,\,.
$$
The cases $k$ odd and $k$ even are slightly different ; we will treat here only 
the case $k$ odd, which is a little more difficult. 
We write $k = 2K + 1$ and decompose $t^{k/2}$ as
\begin{equation*}
\begin{split}
t^{k/2} & = (t-s+s)^K(\sqrt{t} - \sqrt{s} + \sqrt{s}) \\
& = \sum_{l=0}^K \left( \begin{array}{l} K \\ l \end{array} \right)\left[ s^l (t-s)^{K-l} (\sqrt{t} - \sqrt{s}) +  s^l \sqrt{s} (t-s)^{K-l} \right] \,\,,
\end{split}
\end{equation*}
so that we are now dealing with
\begin{equation*}
\begin{split}
t^{k/2} \nabla^k B_2^1 & =
\sum_{l=0}^K \left( \begin{array}{l} K \\ l \end{array} \right) \frac{{\mathbb{P}} \nabla}{\sqrt{-\Delta}}  \int_{0}^{t} (t-s)^{K-l} (\sqrt{t} - \sqrt{s}) \Delta \nabla^{2K-2l+1} e^{(t-s)\Delta} \\
& \;\;\;\;\;\;\;\;\;\;\;\;\;\;\;\;\;\;\;\; \frac{Id-e^{s\Delta}}{\sqrt{-\Delta}} s^l \nabla^{2l} (\chi_{\sqrt{R},x_0} u(s) \times v(s)) \; ds \\
& + \sum_{l=0}^K \left( \begin{array}{l} K \\ l \end{array} \right) \frac{{\mathbb{P}} \nabla}{\sqrt{-\Delta}}  \int_{0}^{t} (t-s)^{K-l} \Delta \nabla^{2K-2l} e^{(t-s)\Delta} \\
& \;\;\;\;\;\;\;\;\;\;\;\;\;\;\;\;\;\;\;\; \frac{Id-e^{s\Delta}}{\sqrt{-\Delta}} s^l \sqrt{s} \nabla^{2l+1} (\chi_{\sqrt{R},x_0} u(s) \times v(s)) \; ds \\
\end{split}
\end{equation*}
For sake of simplicity, let us drop the operator $\frac{{\mathbb{P}} \nabla}{\sqrt{-\Delta}}$, 
which is bounded on $L^2$, and will not play any role, and let us denote $M(s)$ the tensor product 
 $(\chi_{\sqrt{R},x_0} u(s) \times v(s))$. 
The above expression reduces to
\begin{align*}
t^{k/2} \nabla^k B_2^1 
& = \sum_{l=0}^{K-1} \left( \begin{array}{l} K \\ l \end{array} \right) 
Q_{K-l} \left( \frac{Id-e^{s\Delta}}{\sqrt{-\Delta}} s^l \nabla^{2l} M(s) \right) \\
& + Q_{0} \left( \frac{Id-e^{s\Delta}}{\sqrt{-\Delta}} s^K \nabla^{2K} M(s) \right) \\
& + \sum_{l=0}^{K-1} \left( \begin{array}{l} K \\ l \end{array} \right) 
P_{K-l} \left( \frac{Id-e^{s\Delta}}{\sqrt{-\Delta}} s^l\sqrt{s} \nabla^{2l+1} M(s) \right) \\
& + P_{0} \left( \frac{Id-e^{s\Delta}}{\sqrt{-\Delta}} s^K \sqrt{s} \nabla^{2K+1} M(s) \right)\,\,.
\end{align*}

Here we shall concentrate on two types of terms: 
\begin{enumerate} 
\item[(a)] $P_{K-l} \left( \frac{Id-e^{s\Delta}}{\sqrt{-\Delta}} s^l \sqrt{s} \nabla^{2l+1} M(s) \right)$ with $l$ fixed and such that $l < K$. 
\item [(b)] $P_{0} \left( \frac{Id-e^{s\Delta}}{\sqrt{-\Delta}} s^K \sqrt{s} \nabla^{2K+1} M(s) \right)$
\end{enumerate} 
The other terms are dealt with in a very similar way thanks to Proposition \ref{laproposition}.

First, let us concentrate on a term of the form (a),
$$P_{K-l} \left( \frac{Id-e^{s\Delta}}{\sqrt{-\Delta}} s^l \sqrt{s} \nabla^{2l+1} M(s) \right);$$ Our aim is to estimate the $L^2$ norm of this term on a parabolic cylinder. 
Using the boundedness of $P_{K-l}$ on $L^2((0,T), L^2({\mathbb{R}}^d))$ 
for every $T \in (0, \infty]$ (Proposition~\ref{laproposition}), we have
\begin{equation*}
\begin{split}
\frac{1}{R^{\frac{d}{2}}} \int_{0}^{R} \int_{B(x_0,\sqrt{R})} 
& \left| P_{K-l} \left( \frac{Id-e^{s\Delta}}{\sqrt{-\Delta}} 
s^l \sqrt{s} \nabla^{2l+1} M(s) \right) \right|^2 \,dx\,ds \\
& \leq p(K-l) \frac{1}{R^{\frac{d}{2}}} \int_{0}^{R} \int_{{\mathbb{R}}^d} 
\left| \frac{Id-e^{s\Delta}}{\sqrt{-\Delta}} s^l \sqrt{s} \nabla^{2l+1} M(s) \right| ^2 \, dx \,ds  \,\,.
\end{split}
\end{equation*}
Since $\frac{1-e^{-s|\xi|^2}}{|\xi|}$ is bounded by $C\sqrt{s}$, 
the norm of $\frac{(Id-e^{s\Delta})}{\sqrt{-\Delta}}$ as an operator on $L^2$ 
is bounded by $C\sqrt{s}$. Therefore the above term can be bounded by
$$
p(K-l) \frac{1}{R^{\frac{d}{2}}} \int_{0}^{R} \int_{{\mathbb{R}}^d} 
\left| s^{l+1} \nabla^{2l+1} M(s) \right| ^2 \, dx \,ds \,\,.
$$
Observing that $s^{l+1} \nabla^{2l+1} M(s)$ is a sum of terms of the type 
\begin{equation}
\label{MMM}
s^{\frac{m+1}{2}} \nabla^m u\;\;\, s^{\frac{\rho}{2}} \nabla^{\rho}v \;\;\,   
s^{l+1/2-m/2-{\rho}/2} \nabla^{2l+1-m- {\rho}} \chi_{\sqrt{R},x_0} \,\,,
\end{equation}
with $m, \rho \geq 0$ and $m+\rho \leq 2l+1$, we see that
$$
p(K-l) \frac{1}{R^{\frac{d}{2}}} \int_{0}^{R} \int_{{\mathbb{R}}^d} 
\left| s^{l+1} \nabla^{2l} M(s) \right|^2 \; dx \; ds 
\lesssim p(K-l) \sum_{m+\rho \leq 2l+1} \|u\|_{N^m_\infty}^2 \|v\|_{N^{\rho}_C}^2 \,\,.
$$

Now let us consider the term of the type (b), 
$$P_{0} \left( \frac{Id-e^{s\Delta}}{\sqrt{-\Delta}} s^K \sqrt{s} \nabla^{2K+1} M(s) \right).$$
In a similar way to the argument above, we can use Proposition \ref{laproposition}, to obtain: 
\begin{align*} 
& \frac{1}{R^{\frac{d}{2}}} \int_{0}^{R} \int_{B(x_0,\sqrt{R})} \left| 
P_{0} \left( \frac{Id-e^{s\Delta}}{\sqrt{-\Delta}} s^K \sqrt{s} \nabla^{2K+1} M(s) \right)
\right|^2 \,dx\,ds \\
& \leq p(0)
\frac{1}{R^{\frac{d}{2}}} \int_{0}^{R} \int_{B(x_0,\sqrt{R})} \left| 
\frac{Id-e^{s\Delta}}{\sqrt{-\Delta}} s^K \sqrt{s} \nabla^{2K+1} M(s) \right|^2 \,dx\,ds \\
& \leq p(0) 
\frac{1}{R^{\frac{d}{2}}} \int_{0}^{R} \int_{B(x_0,\sqrt{R})} \left| 
 s^{K+1} \nabla^{2K+1} M(s) \right|^2 \,dx\,ds \\
& \leq p(0) 
\left( \|u\|_{N^{0}_{\infty}} \|v\|_{N^{2K+1}_C} + 
 \|v\|_{N^{0}_{\infty}} \|u\|_{N^{2K+1}_C} \right) + 
r(K) \|u\|_{\widetilde{X}^{2K}} \|v\|_{\widetilde{X}^{2K}}
\end{align*} 
We emphasize that in this case the constant $p(0)$ does not depend on $K$. 

Going back to the $L^2$ norm of $t^{k/2} \nabla^k B^1_2$ on a parabolic cylinder, 
and examining precisely the constants on the right-hand side, 
we see that the above inequalities yield
\begin{equation*}
\begin{split}
\frac{1}{R^{\frac{d}{2}}} 
\int_{0}^{R} 
& \int_{B(x_0,\sqrt{R})} |t^{\frac{k}{2}} \nabla^k B_2^1(u,v)(x,t)|^{2} \; dx \; dt \\
& \;\;\;\;\;\;\;\;\leq p_1 \|u\|_{X^0} \|v\|_{\widetilde{X}^k} + p_1 \|v\|_{X^0} \|v\|_{\widetilde{X}^k} + p(k) \|u\|_{\widetilde{X}^{k-1}} \|v\|_{\widetilde{X}^{k-1}} \,\,,
\end{split}
\end{equation*}
where $p_1$ does not depend on $k$; this is the desired bound.

\subsubsection{Estimate for $B^2_2$}

Now we shall prove the following estimate for $B_{2}^{2}$
\be{desB3}
\frac{1}{R^{\frac{d}{2}}}
\int_{0}^{R} \int_{B(x_0,\sqrt{R})} |t^{\frac{k}{2}} \nabla^k B_{2}^{2}|^2 \; dx \; dt 
\lesssim b(k) \; \|u\|_{X^0}^2 \|v\|_{X^0}^2. 
\end{equation} 

We start by applying the boundedness of the Riesz transform on $L^{2}_{t} L^{2}_{x}$ to 
obtain 
\be{B3Ries}
\int_{0}^{R} \int_{B(x_0,\sqrt{R})} |t^{\frac{k}{2}} \nabla^k B_{2}^{2}|^2 \; dx \; dt 
\lesssim 
\int_{0}^{R} \int_{{\mathbb{R}}^{d}} |t^{\frac{k}{2}} (-\Delta)^{\frac{k+1}{2}} 
e^{t \Delta} \int_{0}^{t} M(x,s) \; ds |^{2} \; dx \; dt,
\end{equation}  
where $M(x,s)$ denotes $\chi_{\sqrt{R}, x_0} u(x,s) \times v(k,s)$.

We perform the following change of variables:
$$t = R \tau, \; \; s = R \theta, \; \; x = \sqrt{R} \; z$$ 
to obtain 
\begin{align}
& \int_{0}^{R} \int_{{\mathbb{R}}^{d}} |t^{\frac{k}{2}} (-\Delta)^{\frac{k+1}{2}} 
e^{t \Delta} \int_{0}^{t} M(x,s) \; ds |^{2} \; dx \; dt \nonumber \\
& = \int_{0}^{1} \int_{{\mathbb{R}}^{d}}
| (R \tau)^{\frac{k}{2}} (-\frac{1}{R} {\Delta}_{z})^{\frac{k+1}{2}} e^{\tau {\Delta}_{z}} 
\int_{0}^{\tau} M(\sqrt{R} z, R \theta) \; R d\theta|^2
(\sqrt{R})^d dz \; R d\tau \nonumber \\ 
& =  R^{2 + d/2} \int_{0}^{1} \int_{{\mathbb{R}}^{d}} 
| \tau^{\frac{k}{2}} (-\Delta)^{\frac{k+1}{2}} e^{\tau \Delta} 
\int_{0}^{\tau} N(z, \theta) d\theta |^2
dz \; d\tau \label{cvar1}
\end{align}
where 
$$N(z,\theta) = M(\sqrt{R} z, R \theta).$$ 

In order to obtain an upper bound on the integral in \eqref{cvar1} 
we apply the Lemma \ref{mLe16.2}. Hence 
\begin{equation} \label{mL16.2} 
\int_{0}^{1} \int_{{\mathbb{R}}^{d}} 
|\tau^{\frac{k}{2}} (-\Delta)^{\frac{k+1}{2}} e^{\tau \Delta} 
\int_{0}^{\tau} N(z, \theta) d\theta |^2
dz \; d\tau
\lesssim b(k) A(N) \; \int_{0}^{1} \int_{{\mathbb{R}}^{d}}
 |N(z, \tau)| \; dz \; d\tau,
\end{equation}
where 
$$A(N) =\sup_{x_0 \in {\mathbb{R}}^{d}} \sup_{0 < \tau < 1} \tau^{-\frac{d}{2}}
\int_{0}^{\tau} \int_{|z-x_0| < \sqrt{\tau}} |N(z,\theta)| \; dz \; d\theta.$$

By combining \eqref{B3Ries}, \eqref{cvar1} and \eqref{mL16.2} we have 
\begin{equation} \label{towB3}  
\frac{1}{R^{\frac{d}{2}}}
\int_{0}^{R} \int_{B(x_0,\sqrt{R})} |t^{\frac{k}{2}} \nabla^k B_{2}^{2}|^2 \; dx \; dt
\lesssim R^2  \; b(k) \; A(N) \; \int_{0}^{1} \int_{{\mathbb{R}}^{d}}
 |N(z, \tau)| \; dz \; d\tau. 
\end{equation} 

We proceed by obtaining an upper bound on $\|N\|_{ L^{1}((0,1) \times {\mathbb{R}}^{d}) }$.
After we perform the change of variables 
$$ \sqrt{R} z = x, \; \; R \tau = t$$
we have 
\begin{align} 
& \|N\|_{ L^{1}((0,1) \times {\mathbb{R}}^{d}) } \nonumber \\ 
& = \int_{0}^{1} \int_{{\mathbb{R}}^{d}}
|M(\sqrt{R}z, R\tau)| \; dz \; d\tau \nonumber \\
& = \frac{1}{ R^{1+\frac{d}{2}} } \int_{0}^{R} \int_{{\mathbb{R}}^{d}}
|M(x,t)| \; dx \; dt \nonumber \\
& = \frac{1}{ R^{1+\frac{d}{2}} } \int_{0}^{R} \int_{{\mathbb{R}}^{d}}
\chi_{\sqrt{R},x_0}(x) |u(x,t)| |v(x,t)| \; dx \; dt \nonumber \\
& \leq R^{-1}  
( \frac{1}{R^{\frac{d}{2}}} \int_{0}^{R} \int_{B(x_0,\sqrt{R})} 
|u|^{2} \; dx\; dt)^{\frac{1}{2}} 
( \frac{1}{R^{\frac{d}{2}}} \int_{0}^{R} \int_{B(x_0,\sqrt{R})} 
|v|^{2} \; dx \; dt)^{\frac{1}{2}} \nonumber \\
& = R^{-1} \; \|u\|_{X^0} \|v\|_{X^0}. \label{Nbound}  
\end{align}

In order to obtain an upper bound on $A(N)$ we perform the change of variables
$$  \sqrt{R} z = x, \; \; R \theta = s$$
to obtain 
\begin{align} 
A(N) & =\sup_{x_0 \in {\mathbb{R}}^{d}} \sup_{0 < \tau < 1} \tau^{-\frac{d}{2}}
\int_{0}^{\tau} \int_{|z-x_0| < \sqrt{\tau}} |N(z,\theta)| \; dz \; d\theta  \nonumber \\
& = \sup_{x_0 \in {\mathbb{R}}^{d}} \sup_{0 < \tau < 1} \tau^{-\frac{d}{2}}
\frac{1}{R^{1+\frac{d}{2}}}
\int_{0}^{R\tau} \int_{|\frac{x}{\sqrt{R}}-x_0| < \sqrt{\tau}} |N(\frac{x}{\sqrt{R}},\frac{s}{R})| \; dx \; ds \nonumber \\
& = \sup_{{\widetilde{x_0}} \in {\mathbb{R}}^{d}} \sup_{0 < \tau < 1} \tau^{-\frac{d}{2}}
\frac{1}{R^{1+\frac{d}{2}}}
\int_{0}^{R\tau} \int_{|x-{\widetilde{x_0}}| < \sqrt{R\tau}}
|M(x,s)| \; dx \; ds \nonumber \\
& = \sup_{{\widetilde{x_0}} \in {\mathbb{R}}^{d}} \sup_{0 < \tau < 1} \tau^{-\frac{d}{2}}
\frac{1}{R^{1+\frac{d}{2}}}
\int_{0}^{R\tau} \int_{|x-{\widetilde{x_0}}| < \sqrt{R\tau}}
\chi_{R,x_0}(x) |u(x,s)| |v(x,s)| \; dx \; ds \nonumber \\
& = \sup_{{\widetilde{x_0}} \in {\mathbb{R}}^{d}} \sup_{0 < \tau < 1}
\frac{1}{R} 
( \frac{1}{(R\tau)^{\frac{d}{2}}} \int_{0}^{R\tau} \int_{B({\widetilde{x_0}},\sqrt{R\tau})} 
|u|^{2} \; dx\; dt)^{\frac{1}{2}} 
( \frac{1}{(R\tau)^{\frac{d}{2}}} \int_{0}^{R\tau} \int_{B({\widetilde{x_0}},\sqrt{R\tau})} 
|v|^{2} \; dx \; dt)^{\frac{1}{2}} \nonumber \\
& = R^{-1} \; \|u\|_{X^0} \|v\|_{X^0}. \label{Abound} 
\end{align}

Now we combine \eqref{towB3}, \eqref{Nbound} and \eqref{Abound} to obtain \eqref{desB3}.

\subsection{Conclusion of the argument : proof of Theorem~\ref{reg}}

Proposition~\ref{estimb} gives the following estimate for any $k \geq 1$
\begin{align}
\begin{split} \label{rparticular}
\|B(u,v)\|_{X^k} 
& \leq C_0(k) \|u\|_{X^0} \|v\|_{X^0} \\ 
& + C_1 \|u\|_{X^0} \|v\|_{X^k} + C_1 \|v\|_{X^0} \|u\|_{X^k} \\
& + C(k) \|u\|_{\widetilde{X}^{k-1}} \|v\|_{\widetilde{X}^{k-1}} \,,
\end{split} 
\end{align}
where the constant $C_1$ does not depend on $k$.
On the other hand, the Koch-Tataru solution 
\cite{KT01} satisfies  the following estimate when $k=0$:
\begin{equation} 
\label{rparticular0} 
\|B(u,v)\|_{X^0} \leq C \|u\|_{X^0} \|v\|_{X^0}. 
\end{equation}

Let us define the approximating sequence $u^j$ by 
\begin{align} 
\begin{split} \label{iter} 
& u^{-1} = 0 \\
& u^{0} = e^{t \Delta} u_0 \\
& u^{j+1} = u^0 + B(u^{j},u^j).
\end{split} 
\end{align} 
The usual fixed point argument gives that the sequence $\{u^j\}$ converges in 
$\widetilde{X}^k$ provided that $u^0$ is small enough in $\widetilde{X}^k$. 
But the particular form of the estimates \eqref{rparticular} gives more: 
in the following lemma, we prove convergence of $\{u^j\}$ in $\widetilde{X}^k$ 
for each $k$ under the single condition that $\|u_0\|_{BMO^{-1}}$ is small enough.

\begin{lem}
Let $u_0$ be small enough in $BMO^{-1}$. Then for any $k\geq 0$ there exist constants 
$D_k$ and $E_k$ such that
\begin{equation} \label{estimj}
\| u^j \|_{\widetilde{X}^k} \leq D_k,
\end{equation} 
and 
\begin{equation} \label{estimdifj}  
\|u^{j+1} - u^j\|_{\widetilde{X}^k} \leq E_k \left(\frac{2}{3} \right)^j.
\end{equation}
In particular, for any $k\geq 0$, $u^j$ converges in $\widetilde{X}^k$.
\end{lem}

\begin{proof}
We prove this lemma by induction. When $k=0$ the claim follows from the estimate 
\eqref{rparticular0} via a contraction principle. Also using the estimate 
\eqref{rparticular}, the case $k=1$ follows from a classical contraction argument. 
Furthermore, by choosing $u_0$ small enough in $BMO^{-1}$, we can ensure that for any 
$j \geq 0$ 
\begin{equation}
\label{estimc0}
\| u^j \|_{X^0} \leq \frac{1}{4C_1} \,\,,
\end{equation}
where $C_1$ is the constant in~(\ref{particular}). 

First, let us prove \eqref{estimj}. Assume that \eqref{estimj} is true for $k-1$. 
We shall prove that \eqref{estimj} is true for $k$. 
Applying the estimate~(\ref{rparticular}) to the equation~(\ref{iter}), 
we get
\begin{align} 
\|u^j\|_{X^k} 
& \leq \|u^0\|_{X^k} + \|B(u^{j-1}, u^{j-1} \|_{X^k} \nonumber \\
& \leq \|u^0\|_{X^k} + 2C_1 \|u^{j-1}\|_{X^0} \|u^{j-1}\|_{X^k} 
        +G_k\|u^{j-1}\|^{2}_{\widetilde{X}^{k-1}} \nonumber \\
& \leq \|u^0\|_{X^k} + \frac{1}{2} \|u^{j-1}\|_{X^k} + G_k D_{k-1}^2 
\label{estjind} \\
& \leq D_k + \frac{1}{2} \|u^{j-1}\|_{\widetilde{X}^k}, \label{D}
\end{align} 
where to obtain \eqref{estjind} we used \eqref{estimc0}, the induction hypothesis,
and the notation $G_k = C_0(k) + C(k)$,        
while to obtain \eqref{D} we notice that thanks to Proposition \ref{linear} 
there exists a constant $D_k$ such that 
$\|u^0\|_{X^k} + G_k D_{k-1}^2 \leq D_k$. Hence 
\begin{align*} 
\| u^j \|_{X^k} 
& \leq D_k \sum_{l=0}^{j-1} \left( \frac{1}{2} \right)^l 
+ \frac{1}{2^j} \|u^0\|_{X^k} \\
& \leq D_k \sum_{l=0}^{\infty} \left( \frac{1}{2} \right)^l 
+ \|u^0\|_{X^k} \,\,,
\end{align*}
which is the desired uniform bound (in $j$) on $\| u^j \|_{X^k}$. Thus \eqref{estimj} is proved. 

Now let us prove \eqref{estimdifj}. Considering now the difference of $u^{j+1}$ 
and $u^j$, we have, using the estimate~(\ref{rparticular})
\begin{align} 
\| u^{j+1} - u^j \|_{X^k} 
& \leq \| B(u^j - u^{j-1},u^j ) \|_{X^k} 
+ \| B(u^{j-1} , u^j - u^{j-1} ) \|_{X^k} \nonumber \\
& \leq 
C_1 \|u^j - u^{j-1}\|_{X^k} 
\left( \|u^j\|_{X^0} + \|u^{j-1}\|_{X^0} \right) \label{tr1} \\ 
& + C_1 \|u^j - u^{j-1}\|_{X^0} 
\left( \|u^j\|_{X^k} + \|u^{j-1}\|_{X^k} \right) \label{tr2} \\
& + G_k \|u^j - u^{j-1}\|_{\widetilde{X}^{k-1}} 
\left( \|u^j\|_{\widetilde{X}^{k-1}} + \|u^{j-1}\|_{\widetilde{X}^{k-1}} \right), \label{tr3} 
\end{align} 
where as above $G_k = C_0(k) + C(k)$.
However thanks to \eqref{estimc0} we have 
$$\eqref{tr1} \leq  \frac{1}{2} \|u^j - u^{j-1}\|_{X^k},$$ 
while by applying the induction hypothesis at the level $0$ and \eqref{estimj} we obtain 
$$ \eqref{tr2} \leq  C_1 E_0 \left( \frac{2}{3} \right)^{j-1} 2 D_k  = 3 C_1 E_0 D_k  
 \left( \frac{2}{3} \right)^{j}.$$
On the other hand by the induction hypothesis and \eqref{estimj} 
$$ \eqref{tr3} \leq G_k E_{k-1} \left (\frac{2}{3} \right)^{j-1} 2 D_{k-1} = 
3 G_k E_{k-1} D_{k-1} \left( \frac{2}{3} \right)^{j}.$$
Now let us choose a constant $F_k$ such that 
$$ 3 C_1 E_0 D_k + 3 G_k E_{k-1} D_{k-1} < F_k \,\, .$$ 
Thus we obtain
\begin{align*} 
\| u^{j+1} - u^j \|_{X^k} 
& \leq \frac{1}{2} \|u^j - u^{j-1} \|_{X^k} + F_k \left(\frac{2}{3}\right)^j \\
& \leq \frac{1}{2^{j}} \|u^{1} - u^{0} \|_{X^k} 
+ F_k \sum_{l = 0}^{j-1} \left( \frac{2}{3} \right)^{j-l} \left( \frac{1}{2} \right)^{l} \\
& \leq \left( \frac{2}{3} \right)^j \|u^{1} - u^{0} \|_{X^k} 
+ F_k \left( \frac{2}{3} \right)^{j} \sum_{l = 0}^j \left( \frac{3}{4} \right)^l \\ 
& \leq E_k \left( \frac{2}{3} \right)^j,
\end{align*} 
which proves \eqref{estimdifj} at the rank $k$.
\end{proof}

The lemma which has just been proved states that 
$u^j$ converges to $u$ in $\widetilde{X}^k$, 
and that the $u^j$ are bounded. In particular, 
the norm of $u$ in  $\widetilde{X}^k$ is finite; 
this proves Theorem~\ref{reg}.

\subsection{Proof of Corollary~\ref{tspderiv}}

Let $u$ be Koch-Tataru solution constructed in Theorem \ref{reg}. 
Then 
$$\|\nabla^{\beta} u\|_{L^{\infty}} \lesssim \frac{1}{t^{\frac{\beta+1}{2}}} \mbox{ for } \beta\geq 0,$$ 
which combined with the fact that  $L^{\infty} \hookrightarrow BMO$ gives 
$$\|\nabla^\beta u\|_{BMO} \lesssim 
\frac{1}{t^{\frac{\beta +1}{2}}} \mbox{ for } \beta \geq 0.$$
Now let $\beta = k-1$. The above estimate implies
$$\|\nabla^k u\|_{BMO^{-1}} \lesssim 
\frac{1}{t^{\frac{k}{2}}} \mbox{ for } k \geq 1,$$
and the claim of the colollary is proved for $k \geq 1$.

The case $k=0$ was obtained by Auscher, Dubois and Tchamitchian in \cite{ADT04}. 
For the sake of completeness, here we enclose the proof from \cite{ADT04}. 
Thanks to \eqref{intNS}, Proposition \ref{linear} and the fact that 
the projection into div-free $\mathbb{P}$ maps $L^{\infty}$ to $BMO$, 
it suffices to prove that there exists a positive contact $C$ such 
that for each $t \geq 0$ we have 
\be{0bound} 
\| \int_{0}^{t} e^{(t-s)\Delta} \; (u(\cdot,s) \times v(\cdot,s)) \; ds \|_{L^{\infty}} \leq C. 
\end{equation} 
In a similar way as in Subsection \ref{nkinorm}, we consider separately 
$0 < s < \frac{t}{2}$ and $\frac{t}{2} < s < t$. 

If $0 < s < \frac{t}{2}$ we write the convolution with the heat kernel as follows 
\begin{align*}
& \int_{0}^{t/2} |e^{(t-s)\Delta} (u(x,s) \times v(x,s)) | \; ds\\
& \leq C \frac{1}{t^{\frac{d}{2}}} \int_{0}^{t/2} \int_{{\mathbb R}^{d}} 
e^{ -\frac{|x-y|^2}{4t} } |u (y,s)| |v(y,s)| \; dy \; ds\\  
& \leq C \sum_{q \in {\mathbb Z}^{d}} e^{ -\frac{|q|^2}{10} } \;  
\frac{1}{t^{\frac{d}{2}}} \int_{0}^{t/2} \int_{x-y \in \sqrt{t} (q+[0,1]^{d})} 
|u (y,s)| |v(y,s)| \; dy \; ds, 
\end{align*} 
which after summing in $q$ and applying the Cauchy-Schwartz inequality implies
\be{0bound1} 
\int_{0}^{t/2} |e^{(t-s)\Delta} (u(x,s) \times v(x,s)) | \; ds 
\leq C \|u\|_{X^0} \|v\|_{X^0}.
\end{equation}

However if $\frac{t}{2} \leq s < t$ we have 
\begin{align} 
& \|\int_{t/2}^{t} e^{(t-s)\Delta} (u(\cdot,s) \times v(\cdot,s)) \; ds \|_{L^{\infty}} 
\nonumber \\
& \leq C \int_{t/2}^{t} \|u(\cdot, s)\|_{L^{\infty}} \|v(\cdot, s)\|_{L^{\infty}} \; ds \nonumber \\ 
& \leq C \|u\|_{ N^{0}_{\infty} } \|v\|_{ N^{0}_{\infty} }. \label{0bound2} 
\end{align}

Now \eqref{0bound} follows from \eqref{0bound1} and \eqref{0bound2}. Hence
the claim of the corollary is proved for $k=0$ too.

\section{Analyticity of the solution}

In this section we prove Theorem \ref{analyt} which claims that the solution of 
the Navier-Stokes equation obtained in Theorem \ref{reg} is space analytic, 
meaning its Taylor series converges in $L^{\infty}$ norm. 

{\bf{Proof of Theorem~\ref{analyt}}}
In order to prove analyticity it is enough to show an estimate of the form 
\be{Linfkk} 
\|\nabla^k u\|_{L^{\infty}} \lesssim C^k \frac{k^k}{t^{\frac{k+1}{2}}}. 
\end{equation}

We remark that to obtain \eqref{Linfkk} it suffices to prove 
\be{ssuf} 
\|u\|_{N^k_{\infty}} \lesssim C^{k-1} k^{k-1}, \; \; k \geq 1.
\end{equation} 
We proceed by proving \eqref{ssuf}.

\subsection{Dependence on ${\mathbf{k}}$ in the analysis of the linear term} \label{link} 

Since we want to estimate precisely the growth rate in $k$ of $e^{t\Delta} u_0$
- even if our estimate will not be optimal - some computations are necessary.

\begin{prop} \label{linkk} 
There exists a constant $C$ such that for any $k \in \mathbb{N}$ and $u_0 \in BMO^{-1}$ we have
$$
\| e^{t\Delta} u_0 \|_{N^k_{\infty}} \leq \left( C \sqrt{k} \right)^{k+2} t^{-\frac{k+1}{2}} \|u_0\|_{BMO^{-1}}\,\,.
$$
\end{prop}
So in particular, for $k \geq 3$ we have
$$
\| e^{t\Delta} u_0 \|_{N^k_{\infty}} \leq C^{k-2} k^{\frac{2k}{3}} \|u_0\|_{BMO^{-1}} \,\,.
$$
\begin{proof} 
We will need the $P_j$ Littlewood-Paley operators, which are given by
$$
P_j = \psi \left( \frac{D}{2^j} \right) \,\,,
$$
where 
\begin{equation}
\label{defpsi}
\psi \in {\mathcal{S}}\;\;\;\mbox{is such that}\;\;\; \sum_{j\in \mathbb{Z}} \psi \left( \frac{\xi}{2^j} \right) = 1 \;\;\; \mbox{for}\;\;\; \xi \neq 0 \,\,,
\end{equation}
and $\psi$ supported in an annulus centered in zero.

We begin by estimating the $N^k_\infty$ norm. 
It is well-known that the space $BMO^{-1}$ satisfies
\begin{equation}
\label{BMOBmoinsun}
\sup_{j \in \mathbb{Z}} 2^{-j} \| P_j f\|_{L^{\infty}} \leq C \|f\|_{BMO^{-1}},
\end{equation}
(this is the embedding $BMO^{-1} \hookrightarrow \dot{B}^{-1}_{\infty,\infty}$). 

Now pick an integer $N$ and estimate
\begin{align} 
& \| \nabla^k e^{t\Delta} u_0\|_{L^{\infty}} \nonumber \\
& \leq \sum_{j\leq N} \|P_j\left(\nabla^k e^{t\Delta} u_0\right)\|_{L^{\infty}} 
+ \sum_{j>N} \|P_j\left(\nabla^k e^{t\Delta} u_0\right)\|_{L^{\infty}} \nonumber \\
& \leq \sum_{j\leq N} \|P_j\left(\nabla^k e^{t\Delta} u_0\right)\|_{L^{\infty}} 
+ C \sum_{j>N} 2^{-2j} \|P_j\left(\nabla^{k+2} e^{t\Delta} u_0\right)\|_{L^{\infty}}, \label{lp}
\end{align}
where to obtain \eqref{lp} we used the following propery 
of Littlewood-Paley operators, see, for example, \cite{tao}:
$$ \|\nabla P_j f\|_{L^{\infty}} \sim 2^j \|P_j f\|_{L^{\infty}}.$$

Now \eqref{lp} combined with the boundness of the 
operator $\nabla e^{t\Delta}$ in $L^{\infty}$ given by  
\eqref{heatb} implies 
\begin{align} 
& \|\nabla^k e^{t\Delta} u_0\|_{L^{\infty}} \nonumber \\
& \leq  
\sum_{j\leq N} \| \left( \nabla e^{\frac{t}{k} \Delta} \right)^k P_j u_0 \|_{L^{\infty}}   
+ \sum_{j>N} 2^{-2j} \; \| \left( \nabla e^{\frac{t}{k+2} \Delta} \right)^{k+2} P_j u_0\|_{L^{\infty}}
\nonumber \\
& \leq 
\sum_{j\leq N} \left( C \sqrt{\frac{k}{t}} \right)^k \| P_j u_0 \|_{L^{\infty}}   
+ \sum_{j>N} 2^{-2j} \; \left(C \sqrt{\frac{k+2}{t}} \right)^{k+2} \|P_j u_0\|_{L^{\infty}}
\nonumber \\
& \leq \|u_0\|_{BMO^{-1}} 
\left[ \sum_{j \leq N} \left(C\sqrt{\frac{k}{t}}\right)^k 2^j 
+ \sum_{j>N} \left(C\sqrt{\frac{k+2}{t}} \right)^{k+2} 2^{-j} \right] \label{useinfbmo} \\
& \leq C  \|u_0\|_{BMO^{-1}} 
\left[ \left(C\sqrt{\frac{k}{t}}\right)^k 2^N 
+ \left(C\sqrt{\frac{k+2}{t}}\right)^{k+2} 2^{-N} \right], \label{infbmo-1}
\end{align} 
where to obtain \eqref{useinfbmo} we use \eqref{BMOBmoinsun}. 
Now it suffices to choose $2^N \sim \frac{1}{\sqrt{t}}$ in \eqref{infbmo-1}
to see that 
$$ \|\nabla^k e^{t\Delta} u_0\|_{\infty} 
\leq \left( C \sqrt{k} \right)^{k+2} t^{-\frac{k+1}{2}} \|u_0\|_{BMO^{-1}}\,\,,
$$
which concludes the proof of the proposition. 
\end{proof}

\subsection{Dependence on ${\mathbf{k}}$ in the analysis of the solution}

Here we shall prove
\be{suf} 
\|u\|_{N^k_{\infty}} \leq C^{k-1} k^{k-1}, \mbox{ for all integers } k \geq 1. 
\end{equation}
via mathematical induction. 

First we recall a combinatorial result of Kahane \cite{Kah}: 

\begin{lem} \label{combin} 
Let $\delta > \frac{1}{2}$. Then there exists a constant $C=C(\delta) >0$
such that 
$$\sum_{\gamma \leq \alpha } \binom{\alpha}{\gamma} 
|\gamma|^{|\gamma|-\delta} |\alpha - \gamma|^{|\alpha - \gamma| -\delta} 
\leq C |\alpha|^{|\alpha|-\delta}, 
\mbox{ for all } \alpha \in {\mathbb N}^{d}_{0}.$$
\end{lem}  

As in \eqref{intNS} we write the solution of the Navier-Stokes equations as
\be{intNSrev} 
u(x,t) = e^{t \Delta} u(x,0) - B(u,u)(x,t).
\end{equation} 
We notice that \eqref{suf} is true for $k=1$.  
Also \eqref{intNSrev}, Proposition \ref{linkk} and \eqref{BNinf} imply that 
\eqref{suf} is true for $k=2$. 
Fix $k \geq 3$. Now let us assume that \eqref{suf} is true for $0,1,..., k-1$. 
We shall prove that it is true for $k$ too.

First, we shall prove that  
\be{BNk} 
\|B(u,v)\|_{N^k_{\infty}} \lesssim c_0(k) \|u\|_{X^0} \|v\|_{X^0} + 
c_1 \|u\|_{N^0_{\infty}} \|v\|_{N^k_{\infty}} + c_1 \|v\|_{N^0_{\infty}} \|u\|_{N^k_{\infty}} 
+ C^{k-2} k^{k-1}. 
\end{equation} 
In order to prove \eqref{BNk} we revisit \eqref{BNinf} which gives:
\begin{align} 
\begin{split} \label{improv} 
\|B(u,v)\|_{N^k_{\infty}} 
& \lesssim c_0(k) \|u\|_{X^0} \|v\|_{X^0} + 
c_1 \|u\|_{N^0_{\infty}} \|v\|_{N^k_{\infty}} + c_1 \|v\|_{N^0_{\infty}} \|u\|_{N^k_{\infty}} \\
& + \sum_{l=1}^{k-1} \binom{k}{l}\|u\|_{N^l_{\infty}} \|u\|_{N^{k-l}_{\infty}}. 
\end{split} 
\end{align} 
However the assumption of the mathematical induction combined with Lemma \ref{combin} 
gives that 
\begin{align} 
\sum_{l=1}^{k-1} \binom{k}{l} \|u\|_{N^l_{\infty}} \|u\|_{N^{k-l}_{\infty}} 
& \lesssim \sum_{l=1}^{k-1} \binom{k}{l} C^{l-1} l^{l-1} \; C^{k-l-1} (k-l)^{k-l-1} \nonumber \\
& \leq  C^{k-2} k^{k-1}, \label{bycomb}
\end{align} 
which together with \eqref{improv} implies the desired bound \eqref{BNk} on $B(u,v)$.

Now we combine the remark following Proposition \ref{linkk} and \eqref{BNk} to obtain 
\begin{align} 
\|u\|_{N^k_{\infty}} 
& \lesssim C^{k-2} k^{\frac{2k}{3}} \|u_0\|_{BMO^{-1}} +  c_0(k)\|u\|_{X^0} \|u\|_{X^0} 
+ 2c_1 \|u\|_{N^0_{\infty}} \|u\|_{N^k_{\infty}} + C^{k-2} k^{k-1} \nonumber\\
& \lesssim C^{k-2} k^{k-1} \|u_0\|_{BMO^{-1}} +  c_0(k) \|u\|_{X^0} \|u\|_{X^0} 
+ 2c_1 \|u\|_{N^0_{\infty}} \|u\|_{N^k_{\infty}} + C^{k-2} k^{k-1} \label{uNk}, 
\end{align} 
since $\frac{2k}{3} \leq k-1$ for $k \geq 3$. 
However  $\|u\|_{N^0_{\infty}}$ is small so that the term 
$2c_1 \|u\|_{N^0_{\infty}} \|u\|_{N^k_{\infty}}$
can be incorporated into the left hand side
of \eqref{uNk}. Also our choice of $m$ given in 
\eqref{them} implies that $c_0(k)$ (see \eqref{ms-cons})
is like $C^{k-2} k^{k-1}$. Hence \eqref{suf} is proved.

\section{Proof of Theorem~\ref{thselfsim}}

\label{proofselfsim}

Let $u_0$ be self-similar data with a small norm in $BMO^{-1}$. 
By the Koch and Tataru theorem, there exists a unique solution $u$ small 
enough in $X^0$. Using this uniqueness property and the scaling invariance 
of $BMO^{-1}$ and $X^0$, we see that for any $\lambda$
$$
u(x,t) = \lambda u (\lambda^2 t , \lambda x) \,\,.
$$
Since $u$ is a weakly continuous function with values 
in $BMO^{-1}$ (see Dubois~\cite{D02}), it makes sense 
to define $\psi = u(1,\cdot)$ and
the above equality after taking $\lambda = \frac{1}{\sqrt{t}}$ becomes
$$
u(x,t) = \frac{1}{\sqrt{t}} \psi \left( \frac{x}{\sqrt{t}} \right) \,\,.
$$
It only remains to prove the regularity of $\psi$. 
The $L^\infty$ bound in~(\ref{estimselfsim}) is obvious. As for the other one, we have, by
to the definition of $X^k$,
\begin{equation}
\label{orange}
\int_{0}^{R} \int_{B(x_0,\sqrt{R})}  
t^k \left| \nabla^k u(y,t) \right|^{2} \, dy \, dt \leq C R^{d/2} \,\,.
\end{equation}
Let us now replace $u(y,t)$ by $\frac{1}{\sqrt{t}} \psi \left( \frac{y}{\sqrt{t}} \right)$ and $x_0$ by $0$
in the left hand side of the above inequality:
\begin{equation}
\begin{split}
\int_{0}^{R} \int_{B(0,\sqrt{R})}  
t^k \left| \nabla^k u(y,t) \right|^{2} \, dy \, dt 
& = \int_{0}^{R} \int_{B(0,\sqrt{R})} 
\frac{1}{t} \left| \nabla^k \psi \left( \frac{y}{\sqrt{t}}
\right) \right|^{2} \, dy \, dt \\
& = \int_{0}^{R} \int_{B(0,\sqrt{R/t})} 
t^{d/2-1}  \left| \nabla^k \psi(y) \right|^{2} \, dx \, dt \\
& = \int_{{\mathbb{R}}^d} \left| \frac{\sqrt{R}}{y} \right|^d 
\left| \nabla^k \psi \right|^{2} \, dy\;. \label{calcself}
\end{split}
\end{equation}
Now we combine \eqref{orange} with \eqref{calcself} to obtain 
$$
\int_{{\mathbb{R}}^d} \left| \frac{\sqrt{R}}{y} \right|^d 
\left| \nabla^k \psi \right|^{2} \, dy \leq C R^{d/2} \,\,,
$$
which is the desired result.

\end{document}